\input amssym.def
\input amssym.tex
 

\magnification=\magstep1

\hfuzz=10pt
\font\bigrm=cmb10 scaled 1200
\font\nrm=cmcsc10 at10pt
\font\prm=cmr10 at9pt
\font\pit=cmsl10 at9pt
\font\pbo=cmb10 at9pt
\font\psl=cmsl10 at9pt

\def\ul{\underline}

\def\ra{\rightarrow}
\def\lra{\longrightarrow}
\def\geq{\geqslant}
\def\leq{\leqslant}

\def\s{{\sigma}}
\def\om{{\omega}}

\def\im{{\frak m}}
\def\ib{{\frak b}}

\def\proj{{\rm Proj}}
\def\codim{{\rm codim}}

\def\reg{{\rm reg}}
\def\ann{{\rm Ann}}

\def\tor{{\rm Tor}}

\def\Z{{\cal Z}}
\def\Zi{{{\cal Z}_{I}}}
\def\Sc{{\cal S}}
\ \bigskip\bigskip
\centerline{\bigrm Projective schemes: What is Computable in low
            degree?}
\bigskip
\centerline{Marc Chardin}

\bigskip\bigskip\bigskip

{\it To Professor Wolmer Vasconcelos, who inspired this work in many ways.}

\bigskip\bigskip
\centerline{\bf Introduction}\bigskip

Let $I$ be a homogeneous ideal in $A:=k[X_{0},\ldots ,X_{n}]$ (where
$k$ is a field) given in terms of its generators
$$
I=(f_{1},\ldots ,f_{s})
$$
where $f_{i}$ is a form of degree $d_{i}$.

This $I$ defines a scheme
$$
\Zi \subset {\bf P}_{n}(k),
$$
and there is a one to one correspondance between the subschemes of ${\bf
P}_{n}(k)$ and the homogeneous ideals up to saturation (the
saturation $I^{*}$ of $I$ consists of elements $f$ in $A$ such that
for some $m$, $X_{i}^{m}f$ is in $I$ for all $i$). 

We will mention few ideas for algorithms to compute geometric informations
on $\Zi $ from the generators of $I$ and speak about one aspect of their
complexity. Algebraic geometry told us that many geometric invariants
may be computed from objects that have a more algebraic flavor: finite
free resolutions, cohomology, Hilbert function, etc. 

There are many ways of splitting the algorithmic problems into parts,
we will choose the following one:\smallskip

1) Provide algorithms that are easy to program,\smallskip
2) Estimate their complexity in terms of the output and/or the
input,\smallskip
3) Bound the complexity of the output in terms of the input.\smallskip

As we are not at all expert in complexity theory, we will choose a
measure of complexity that we know : Castelnuovo-Mumford 
regularity. It bounds the degree (in $A$) where most algebraic
questions reduces to linear algebra problems. 
\bigskip

{\bf 1. Main ingredients of two simple algorithms.}\medskip

We consider $I=(f_{1},\ldots ,f_{s})$ a homogeneous ideal in
$A:=k[X_{0},\ldots ,X_{n}]$ ($k$, a field) set $d_{i}:=\deg f_{i}$ and
assume for simplicity that $d_{1}\geq \cdots \geq d_{s}\geq
1$ and that $k$ is infinite.\medskip

We choose two very simple algorithms as illustrations of what we look
for, there are more details (and other algorithms) in [Ch1] for the
first one and the second is based on a particular case of [Ch2, 5.2]. 

The first one relies in part on the following lemma ([Ch1,
20])\medskip

{\bf Lemma 1.1.} {\sl Let $g_{1},\ldots ,g_{t}$ be a homogeneous regular
sequence in $A$ with $g_{i}=f_{i}+\sum_{j>i}h_{i,j}f_{j}$ and let $J$
be the ideal they generate. Set  $\s :=d_{1}+\cdots +d_{r}-r$. The
following  are equivalent, \smallskip  

{\rm (1)} $\codim (I)>t$,\smallskip 

{\rm (2)} there exits $g_{t+1}=f_{t+1}+\sum_{j>t+1}h_{t+1,j}f_{j}$ such that 
 $g_{1},\ldots ,g_{t+1}$ is an homogeneous regular
sequence in $A$,\smallskip 

{\rm (3)} the map
$$
(A/J)_{\s}\buildrel{(f_{t+1},\ldots ,f_{s})}\over{\lra}
\bigoplus_{j=t+1}^{s}(A/J)_{\s +d_{j}}
$$
is injective.}
\medskip

{\bf Steps of algorithm 1:}

- Step 1: Construct a sequence $g_{1},\ldots ,g_{r}$ as in the lemma
  with $r=\codim (I)$, using (3) to determine if $t=r$ and elementary
  transformations in the matrix representing this $k$-linear map to
  construct $g_{t+1}$ as in (2) if $t<r$. 

- Step 2: Choose an homogeneous element $h$ in the kernel of  
$$
(A/J)\buildrel{(f_{r+1},\ldots ,f_{s})}\over{\lra}
\bigoplus_{j=r+1}^{s}(A/J)[d_{j}]
$$
of degree at most $\s :=d_{1}+\cdots +d_{r}-r$ (there exists such an
element by (3) and the minimal degree of such an element provides an
interesting invariant of $\Zi $: the $a$-invariant of the coordinate
ring of $\Zi $).  

- Step 3: Compute the kernel 
$$
(A/J)\buildrel{\times h}\over{\lra}
(A/J)[\deg h].
$$

The output of the algorithm is the defining ideal of a scheme
$\Sc \subset \Zi $ which is purely of dimension $\dim \Zi $ ({\it i.e.} unmixed
of codimension $r$). If $h$ is ``general'' the support of $\Sc $ is the 
unmixed part of $\Zi $. The complexity is bounded by the following
result (in terms of Gr{\"o}bner basis),\medskip

{\bf Lemma 1.2.}[Ch1] {\sl Set $\ib :=(g_{1},\ldots ,g_{r},h-T^{\deg
  h})\subset A[T]$. A Gr{\"o}bner basis ${\frak B}$ of $\ib$ for the
  deg-rev-lex order ``contains'' Gr{\"o}bner bases for $\ib +(h)$ and
  $I_{\Sc }=\ib :(h)$. The maximal degree of an element in 
${\frak B}$, for general coordinates in the $X_{i}$'s, is at most
$$
\max \{ d_{1}+\cdots +d_{r}-r,\reg(\Sc )+\deg h\} ,
$$
except possibly for the element $T^{2\deg h}$.}
\medskip
In characteristic zero, this bound is in fact achieved for ``very
general'' coordinates.
\medskip
{\bf Steps of algorithm 2:} Assume that $I$ is the defining ideal of normal
scheme $\Sc $ ({\it i.e.} $I=I^{*}$ and $\Sc =\proj (A/I)$).  

- Step 1: Choose two elements $f,g$ in the Jacobian ideal of $I$ that such that
$\codim (I+(f,g))=\codim (I)+2$. 

- Step 2: Compute the $A/I$-module 
$$
H^{1}(f,g;A/I)=\{ (x,y),\ fx+gy=0\} /\{ a(g,-f),
\ a\in A\} .
$$

The output of the algorithm is the local cohomology module $H^{1}_{\im}(A/I)\simeq
H^{1}(f,g;A/I)$ (called the Hartshorne-Rao module). In practice
choosing $f$ and $g$ should be easy, verifying the codimension  
condition costs quite a lot as the degrees of $f$ and $g$ are not that small when
the codimension increases. 

Another strategy may be to compute first the last degree in
which $H^{1}_{\im}(A/I)$ is not zero (apply the same algorithm, replacing $f$ and $g$
by two linear forms satisfying the same codimension condition, the last non zero degree 
is the same) and then use linear algebra to finish the computation (in
place of a Gr{\"o}bner basis computation, that doesn't require {\it a
  priori} bounds).\medskip  

The main common point of these algorithmes is that they produce (at
least in some important cases) a module that is encoding
geometric informations with a complexity controled mainly by the
complexity of the ouput.\medskip

{\bf 2. Castelnuovo-Mumford regularity}\medskip

There are many ways to define this invariant attached to a finitely generated
graded module over a polynomial ring. Let us recall some of them in a proposition and 
then connect it to degrees of element in a Gr{\"o}bner basis.
\medskip
{\bf Definition.} Let $A$ be a polynomial ring over a noetherian ring $k$, $M$ be
a finitely generated $A$-module that is graded (for the standard grading of $A$), and
$\im$ the ideal generated by the variables. For an 
integer $i$ we set
$$
a_{i}(M):=\max \{\mu\ \vert\ H^{i}_{\im}(M)_{\mu}\not= 0\}
$$
and
$$
b_{i}(M):=\max \{\mu\ \vert\ \tor^{A}_{i}(M,k)_{\mu}\not= 0\}
$$
(with the convention $\max \emptyset =-\infty$).\medskip

We recall that $a_{i}(M)$ is indeed finite (Serre's vanishing theorem)
and that the Tor module is equal to the Koszul homology module
$H_{i}({\ul x};M)$, where ${\ul x}$ denote the set  of variables
(because $K_{\bullet}({\ul x};A)$ provides a free $A$-resolution of
$A/\im =k$ as an $A$-module). Also, if $k$ is a field and
$F_{\bullet}\ra M\ra 0$ is a minimal free resolution of $M$, then 
$\tor^{A}_{i}(M,k)_{\mu}:=H_{i}(F_{\bullet}\otimes _{A}A/\im
)=(F_{i}\otimes _{A}A/\im )_{\mu}$ (maps in $F_{\bullet}$ are
represented by matrices with entries in $\im$) is the 
number of minimal  generators of degree $\mu$ of $F_{i}$, because
$A[-j]\otimes _{A}A/\im =k[-j]$ is concentrated in degree $j$.
\medskip
(See [Ei, Ch. 17] or [BH, Ch. 1] for the definition and basic facts on
the Koszul complexes $K_{\bullet}({\ul z};M)$, $K^{\bullet}({\ul
  z};M)$ associated to a module $M$ and a tuple ${\ul z}$ of elements
of $A$. We will denote by $H_{i}({\ul z};M)$ and $H^{i}({\ul
  z};M)$ their homology (resp. cohomology) modules.)

For simplicity, we will assume the base ring to be a field in the following,\medskip

{\bf Proposition 2.1.} {\sl Let $A$ be a polynomial ring over a field $k$ and $M$ be
a finitely generated $A$-module that is graded (for the standard grading of $A$). 

The following definitions are equivalent,\smallskip

{\rm (1)} $\reg (M):=\max_{i}\{ a_{i}(M)+i\}$,\smallskip

{\rm (2)} $\reg (M):=\max_{i}\{ b_{i}(M)-i\}$,\smallskip

{\rm (3)} $\reg (M):=\min\{ \mu \ \vert\ b_{i}(M_{\geq \mu})\leq b_{0}(M_{\geq \mu})+i\ \forall i\}$,\smallskip

{\rm (4)} Let ${\ul z}$ be a finite collection of homogeneous elements of $A_{>0}$ such that $M/({\ul z})M$ is a finite 
dimensional vector space,
$$
\reg (M):=\min\{ \mu \ \vert\ H^{i}({\ul z};M)_{>\mu -i}=0\ \forall i\} .
$$
}

The equivalence of (1) and (4) is easy using the standard tool for
comparing two homological objects (spectral sequences), it implies the
equivalence with (2) (taking the set of variables for ${\ul z}$). Using the
equivalence of (1) and (2) and the fact that $a_{i}(M)=a_{i}(M_{\geq
  \mu})$ for any $\mu$ and $i>0$, leads to the equivalence with (3). 
\medskip

The definition (1), despite its apparent inaccessibility (e.g. the
local cohomology modules are not all finitely  generated, unless $M$
is of finite length) happens to be most tractable when one wants to
estimate the regularity. As always, some familiarity with the object
makes them very concrete and most of their apparent pathologies are
not so bad (e.g. the graded duals of local cohomology modules are
finitely generated, as they  are isomorphic to some Ext
modules). \medskip 

Definition (2) is interesting for looking at
families of schemes (the study of the Hilbert scheme). The condition
is that the maps in the minimal free resolution of $M_{\geq \mu}$ over
$A$ have linear forms as entries. The sheaves associated to $M$
and $M_{\geq \mu}$ are the same. The existence of {\it a priori}
bounds on the regularity is one ingredient for proving the existence
of Hilbert schemes (see [Mu], or [EH, Ch. VI] for an introduction).

Note that definition (4) implies in particular that if $M$ is of
dimension $d$ and is a finitely generated $B$-module where
$B:=k[l_{1},\ldots ,l_{d}]$ and the $l_{i}$'s are linear forms (in
other  words, we have a Noether normalisation) then the regularity of
$M$ as a $B$-module is the same as the one as an $A$-module --this is
also quite immediate from (1).  \medskip 

One way to connect the regularity to degree of generators of a Gr{\"o}bner
basis, is to study how it behaves when passing modulo a ``general''
linear form (see e.g. [Ei, 20.20 and 20.21]).  The key is the
following lemma, where $A$ is a polynomial ring over a field,\medskip

{\bf Lemma 2.2.} {\sl Let $M$ be a finitely generated graded $A$-module and
$l$ a linear form. Set  
$K:=0:_{M}(l)=\{ m\in M\ \vert\ lm=0\}$. Then $\reg (M)\leq \max \{
\reg (K),\reg (M/(l)M)\}$,  
and if the Krull dimension of $A/\ann (K)$ is at most one, then
$$
\reg (M)=\max \{ \reg (K),\reg (M/(l)M)\}.
$$
}

The proof is a standard diagram chasing using the local cohomology
definition of the regularity (which has the advantage that
$H^{i}_{\im}(N)$ for $i>\dim N$). This lemma gives in particular the
following,\medskip 

{\bf Proposition 2.3.} {\sl Let $\Sc =\proj (A/I)$ be a projective surface
  ({\it i.e.} an unmixed scheme of dimension two). Assume that $\Sc
  \cap \{ X_{n}=X_{n-1}=0\}$ is a zero dimensional scheme. Then, for the
  deg-rev-lex order,
$$
\reg (I)=\reg (in(I)).
$$
}

Here $in(I)$ is the ideal generated by the leading monomials of the
polynomials in $I$ for a given order on the monomials. 
We recall that ${\frak B}$ is a (minimal)
Gr{\"o}bner basis of $I$ if $\{ in(f)\ | \ f\in {\frak B}\}$ are (minimal)
generators of $in(I)$. Therefore, the maximal degree of an element in  
a minimal Gr{\"o}bner basis of $I$ is $b_{0}(in(I))$.  The deg-rev-lex 
order on monomials is obtained by refining the degree order by the 
inverse of the lexicographic order. (See e.g. [Ei,
Ch. 15] for an introduction on Gr{\"o}bner bases.)

With no geometric hypotheses on $\proj (A/I)$ one of the main early
discoveries is the following result  of Bayer and Stillman,\medskip 

{\bf Theorem 2.4.}[BS] {\sl For any order and any coordinates, $\reg (I)\leq
\reg (in(I))$. For the deg-rev-lex order, and in general coordinates,
$\reg (I)=\reg (in(I))$.\medskip 
}

The expression ``general coordinates'' means that there exists a
Zariski open subset of the linear group so that any matrix of this
open subset gives rise to coordinates that satisfies the given
property (in  particular it may be that over finite fields an
extension of the base field is needed to find good coordinates). 

Also, Diana Taylor find an explicit resolution of monomial ideals, that
in particular proves the following,\medskip 

{\bf Proposition 2.5.} {\sl If $J$ is a monomial ideal in $A$,
$$
b_{i}(J)\leq (i+1)b_{0}(J),
$$
so that $\reg (J)\leq (n+1)(b_{0}(J)-1)+1$.
}
\medskip
Note that, in this form, the result is optimal: consider the case
where $J$ is generated by the $b_{0}$-th powers of the variables.\medskip

Mayr and Mayer, and other since then, provided examples of binomial
ideals $J$ where $b_{i}(J)$ is much bigger than $b_{0}(J)$ (say $\reg
(J)\gg b_{0}(J)^{c^{n}}$ where $n$ is the number of variables and
$c>1$ is close to $\sqrt{2}$). 

Bounds on the regularity follows from an inductive
argument on the number of variables from the following result (see 
[BM, the proof of 3.8]). The bounds are more or less of the same type
as the lower bounds coming from Mayr-Meyer type examples.  
\medskip 

{\bf Proposition 2.6.} {\sl If $M$ is a finitely generated graded
  $A$-module, for a general linear form $l$ and $\mu \geq \max \{ \reg
  (M/(l)M)+1,b_{0}(M),b_{1}(M)-1\}$, 
$$
\reg (M)\leq \mu +\dim_{k}H^{0}_{\im}(M)_{\mu}\leq \mu
+\dim_{k}M_{\mu}
$$
}

{\it Idea of the proof.} Notice that, for a general $l$ the kernel $K:=\ker
(M\buildrel{\times l}\over{\ra} M)$ is of finite length. Some
diagram chasing gives $b_{0}(K)\leq \max \{ 
b_{0}(M),b_{1}(M)-1,b_{2}(M/(l)M)-1\}$ so that $K$ has no generator of
degree bigger than $N:=\max \{ \reg
(M/(l)M)+1,b_{0}(M),b_{1}(M)-1\}$ (in particular $K_{\mu}=0$ implies
$K_{\mu +1}=0$ for $\mu \geq N$). 

Now the local
cohomology definition shows that $\reg (M/H^{0}_{\im}(M))\leq \reg
(M/(l)M)$ and gives an exact sequence $0\ra K_{\mu}
\ra H^{0}_{\im}(M)_{\mu}\buildrel{\times l}\over{\ra}H^{0}_{\im} (M)_{\mu
+1}\ra 0$ for $\mu > \reg (M/(l)M)$, so that for $\mu \geq N$,
$\dim_{k}H^{0}_{\im}(M)_{\mu}$ is a strictly deacreasing function of
$\mu$ until it reaches 0. 
\medskip

It was also remarked (and proved) by Andr{\'e} Galligo that initial ideals
have an interesting property in general coordinates and characteristic
zero: they are stable, which means that if a monomial $x_{i}m$ is in
the initial then it is also the case of $x_{j}m$ if $j<i$. Afterwards,
Eliahou and Kervaire find a minimal free resolution for stable
monomial ideals; it follows from the resolution that all graded Betti
numbers of these ideals may be easily read from the minimal
generators, in particular $\reg (J)=b_{0}(J)$ for a stable monomial
ideal $J$.  

The conditions of genericity needed for having a stable monomial ideal
are not so often realized without performing a change of coordinates,
they are more diffcult to achieve than the ones for having $\reg (I)=
\reg (in(I))$ (for deg-rev-lex order).\medskip

It is also important to notice that doing a generic change of
coordinates have quit a big influence on the size of the computation,
for several reasons: the coefficients get bigger, the polynomials
became dense and the number of generators of the initial ideal
increases in general. On the other hand it should be noted  that the
degrees of generators of the initial ideal in special coordinates may
be much bigger than in  general coordinates. Applying the following
lemma to the Mayr-Mayer ideal provides such an example,\medskip

{\bf Lemma 2.7.} {\sl Let $f_{1},\ldots ,f_{s}$ be forms of degrees
$d_{1},\ldots d_{s}$ in $A$ and set $I:=(
f_{1},\ldots ,f_{s})$.\smallskip
{\rm (1)} If ${\rm codim}(I)=s$, $\reg (I)=d_{1}+\cdots
+d_{s}-s+1$.\smallskip
{\rm (2)}  If ${\rm codim}(I)=r<s$, there exists a graded complete
intersection of degrees $d_{1},\ldots d_{s}$ in $A[Y_{1},\ldots
,Y_{s-r}]$ such that the regularity of its initial ideal for the
deg-rev-lex order bounds the regularity of the initial ideal of $I$ for
the deg-rev-lex order.\medskip
}
 
The regularity of a graded ideal is always
bounded in terms of the Hilbert function of the ideal, as all the
graded Betti numbers are bounded above by the ones of the 
lex-segment ideal that only depends on the Hilbert function (see
[Bi] and [Hu] for characteristic zero case, [Pa] for the general case,
and [CGP] for a short argument in characteristic zero). If $I$ is
saturated, the regularity of the lex-segment ideal only depends on the
Hilbert polynomial. The regularity of the lex-segment ideal may
be computed ([CM, 1.3 and 2.3]) and leads for example to the
following  bound that is at least as bad as expected...\medskip 

{\bf Corollary 2.8.}[CM]  {\sl Let $I\subset A$ be an ideal generated by 
polynomials of degrees $d_{1},\ldots ,d_{s}$ and let $r$ be the
codimension of $I$. For any admissible monomial order and any
coordinates,   
$$
\reg (in(I))\leq  1+[ d_{1}\cdots d_{s}] ^{2^{n-r}}
$$
if $r\leq n$ and $\reg (in(I))\leq d_{1}+\cdots +d_{n+1}-n$ if
$r=n+1$.
}
\bigskip

{\bf 3. Bounds on Castelnuovo-Mumford regularity}\medskip

There is a famous conjecture that suggests the following bound for
reduced and irreducible schemes:

{\bf Conjecture}[Eisenbud and Goto]. {\sl If $\Sc \subset {\bf P}_{n}$
  is a non degenerate reduced and irreducible scheme, 
$$
\reg (\Sc )\leq \deg \Sc -\codim \Sc.
$$
}
(Non degenerate means $\Sc \not\subset H$ for any hyperplane $H$.)

We recall that if $\Sc :=\proj (A/I)$, $\reg (\Sc ):=\reg
(A/I^{*})=\reg (I^{*})-1$.

This result was known for curves when the conjecture was made. It was
first established for smooth curves by Castelnuovo [Ca], and the for
reduced curves with no degenerate component by Gruson, Lazarsfeld and
Peskine (over a perfect field)in [GLP]. There is some evidence that this may
be true at least for smooth schemes in characteristic zero: it is true
for smooth surfaces (Pinkham  
and Lazarsfeld) and (up to adding small constants) in dimension at
most six, by the work of several people including Lazarsfeld, Ran and
Kwack.  

In any dimension, it was prove by Mumford ([BM]) that in characteristic zero
a smooth scheme $\Sc$ satisfies,
$$
\reg (\Sc )\leq (\dim \Sc +1)(\deg \Sc -1).
$$

In positive characteristic, it follows from theorems that we will
mention below that one has $\reg (\Sc )\leq (\dim \Sc +1)^{2}(\deg \Sc
-1)$, and there are also quite reasonable results for schemes with
isolated singularities. 
\medskip

We now turn to bounds depending on the degrees of generators. As we
mentioned in the preceeding paragraph, there is no reasonable bound
on the regularity without imposing geometric conditions. 
 
Let $I=(f_{1},\ldots ,f_{s})$ be a homogeneous ideal in $A$, 
where $f_{i}$ is a form of degree $d_{i}$. We will assume that
$d_{1}\geq d_{2}\geq \cdots \geq d_{s}\geq 1$. 
Let $\Zi \subset {\bf P}_{n}$ be the scheme defined by $I$ and
$r$ be the codimension of $I$ in $A$, which is also the
one of $\Zi $ as a subscheme of $ {\bf P}_{n}$. Let $\Sc$ be the top 
dimensional  part of $\Zi $ and $Y$ the residual of $\Sc$ in $\Zi $. In
algebraic terms, $I_{\Sc}$ is the intersection of the primary components
of $I$ of codimension $r$, and $I_{Y}:=(I:I_{\Sc})=\{ f\in A\ |\
fI_{\Sc}\subset I_{Y}\}$. 

The first striking result on regularity in these terms is due to
Bertram, Ein and Lazarsfeld:

{\bf Theorem 3.1.}[BEL] {\sl If $\Zi =\Sc$ is smooth of characteristic zero,
$$
\reg (\Sc )\leq d_{1}+\cdots +d_{r}-r,
$$
with equality if and only if $\Sc$ is a complete intersection
of degrees $d_{1},\ldots ,d_{r}$.}
\medskip
This linear bound was generalized in [CU], 
\medskip
{\bf Theorem 3.2.} {\sl Assume that $\Sc$ have at most a one dimensional
  singular locus an is locally a complete intersection outside
  finitely many points. If the residual $Y$ have at most isolated
  singularities and $k$ is of characteristic zero,
$$
\reg (\Sc )\leq d_{1}+\cdots +d_{r}-r,
$$
with equality if and only if $\Zi =\Sc$ is a complete intersection
of degrees $d_{1},\ldots ,d_{r}$.}
\medskip

Note that the defining ideal of $\Sc$ may be computed in low degree by
Algorithm 1, even if the regularity of $I$ is much bigger. 
They rely on liaison theory and use either Kodaira's vanishing theorem
or a result of Karen Smith that enables an induction on the
dimension. 

More recently, we showed several other bounds. They essentially
improve the ones of [CU] in positive characteristic, and provide the
following result:
\medskip
{\bf Theorem 3.3.}[Ch2, 4.4] {\sl Assume that ${\cal Z}$ is an
  isolated component of $\Zi $ that doesn't meet the other components,
  and that $\Zi $ is smooth at all but a finite number of points of
  ${\cal Z}$. Then, 
$$
\reg ({\cal Z})\leq (\dim {\cal Z}+1)(d_{1}+\cdots +d_{r}-r-1)+1.
$$
}

This generalizes the result of [CP] that treats the case where $\dim
\Z =0$. The proof relies on [CP] and a result of Hochster and Huneke,
which implies that the phantom homology (which is, roughly speaking,
the one that vanishes in the Cohen-Macaulay case) is uniformly killed
by the Jacobian ideal. The result then follows by cutting the scheme
${\cal Z}$ by a sequence of parameters in the Jacobian ideal and using
some homological algebra to exploite this uniform vanishing. The
connection between annihilators and vanishig was already remarked and
used to study the so-called $\ell$-Buchsbaum schemes by Miyazaki,
Nagel, Schenzel and Vogel ([Mi], [NS1], [NS2] and [MV]). 
\medskip

Let us also point out the following remark that formalizes
the fact that bounding the regularity in a geometric context is as
difficult as bounding the degree where the Hilbert function becomes a
polynomial, or bounding the degree where every global section is the
restriction of a polynomial.\medskip

{\bf Remark 3.4.}[CM 2.5] Let ${\cal P}$ be a property of embedded projective
schemes and $N(X)$ a numerical invariant attached to such a
scheme $X$. Assume that if $X\subseteq {\bf P}_{n}$ satisfies ${\cal P}$ and
$H$ is a general hyperplane, then $X\cap H\subseteq H\simeq {\bf P}_{n-1}$
satisfies ${\cal P}$ and $N(X\cap H)\leq N(X)$. We denote by
$I_{X}\subseteq R$ the defining ideal of $X$ and by $H_{X}$ the Hilbert
function of $R/I_{X}$.  Then the following are equivalent,\smallskip 

(i) If $X$ satisfies ${\cal P}$, $\reg (R/I_{X})\leq N(X)$.\smallskip

(ii) If $X$ satisfies ${\cal P}$, $\reg (H_{X})\leq N(X)-1$.\smallskip

(iii)  If $X$ satisfies ${\cal P}$, $(R/I_{X})_{\mu}=H^{0}(X,{\cal
O}_{X}(\mu ))$ for $\mu \geq N(X)$,\smallskip

\noindent where $\reg (H_{X})$ is the last degree where $H_{X}$ differs from the
Hilbert polynomial $P_{X}$.\medskip

Examples for property ${\cal P}$ are: $X$ satisfies
$S_{k}$, $X$ is smooth in codimension $\ell$, $X$ is irreducible, $X$ is
equidimensional, or any conjonction of some of these properties. For
$N(X)$ one may choose the degree of $X$, or the 
degree of $X$ minus the embedding codimension of $X$ if $X$ is
irreducible and reduced, or the minimum over the sets of equations defining
$X$ of the maximal degree of these equations.\medskip

Another important point is to notice that even if we are not able to
bound the regularity in many cases, a big part of the information is
sometimes available in an indirect way. For example, if the top
dimensional component $\Sc$ of $\Zi$ have at most isolated
singularities the canonical module $\om_{\Sc}$ of $\Sc$ have a small
regularity (at least in characteristic zero, thanks to Kodaira's
vanishing theorem) and is easily computable (as the kernel of the map
in Step 2 of Algorithm 1). From $\om_{\Sc}$ we may compute the Hilbert
polynomial of $\Sc$ or its cohomology modules using Serre duality (at
least if $\Sc$ is Cohen-Macaulay), or test if an element is in
$I_{\Sc}$. This in turn gives a way to check if $\reg (\Sc )\leq N$ by
linear algebra computation in degree at most $N$ plus a linear
function of the degrees of generators (in general coordinates, the
criterion of [BS] gives such a test).\bigskip

{\bf Acknowledgements.} A major part of my understanding of
Castelnuovo-Mumford  grew out of the many ideas learned from or shared
with several people. It is partially reflected in the references
below: Patrice Philippon and Bernd Ulrich for the ``geometric''
results and Vesselin Gasharov, Irena Peeva and Guillermo Moreno for
our work on lex-segment ideals. Many others should be mentioned: Dale
Cutkosky, David Eisenbud, Philippe Gimenez, Robert Lazarsfeld, Mike
Stillman, Wolmer Vasconcelos, ...    
A special thanks also to Martin Sombra for reading and commenting a
first draft of these notes.

\bigskip\bigskip\bigskip

{\bf References.}\medskip

{\prm 

[BM] D. Bayer, D. Mumford, {\psl 
What can be computed in algebraic geometry?} in 
{\pit Computational algebraic geometry and commutative algebra
  (Cortona, 1991), 1--48}. Sympos. Math. XXXIV, Cambridge
Univ. Press, Cambridge, 1993.\par  

[BS] D. Bayer, M. Stillman, {\psl A criterion for detecting $m$-regularity.}
Invent. Math. {\pbo 87} (1987), 1--11.\par 

[BEL] A. Bertram, L. Ein, R. Lazarsfeld, {\psl
Vanishing theorems, a theorem of Severi, and the equations defining
projective varieties.} J. Amer. Math. Soc. {\pbo 4} (1991), 587--602.\par   

[BH] W. Bruns, J. Herzog. {\pit Cohen-Macaulay Rings.} Cambridge stud. in adv. math
{\pbo 39}. Cambridge Univ. Press, 1993.\par

[Bi] A. Bigatti, {\psl Upper bounds for the Betti numbers of a given
Hilbert function.}
Comm. in Algebra {\pbo 21} (1993), 2317--2334.\par

[Ca] G. Castelnuovo, {\psl Sui multipli di una serie lineare di gruppi di punti
appartenente ad una curva algebrica.}  Rend. Circ. Mat. Palermo
{\pbo 7} (1893), 89--110.\par

[Ch1]  M. Chardin, {\psl
Applications of some properties of the canonical module in
computational projective algebraic geometry}. J. Symbolic Comput. {\pbo
29} (2000), 527--544.\par 

[Ch2]  M. Chardin, {\psl Cohomology of projective schemes: from
  annihilators to vanishing.}({\pit preprint} 290, U. Paris 6)\par

[CGP] M. Chardin, V. Gasharov, I. Peeva, {\psl Maximal Betti numbers}. 
Proc. Amer. Math. Soc. (to appear).\par

[CM] M. Chardin, G. Moreno, {\psl Regularity of lex-segment ideal: some
      closed formulas and applications}, ({\pit preprint} 292, U. Paris 6).\par

[CP] M. Chardin, P. Philippon, {\psl R{\'e}gularit{\'e} et interpolation}.
J. Algebraic Geom. {\pbo 8} (1999), 471--481. \par
{\pit See also the {\rm erratum} at the address} 
http://www.math.jussieu.fr/$\tilde{\;}$chardin/textes.html.\par

[CU] M. Chardin, B. Ulrich, {\psl Liaison and Castelnuovo-Mumford regularity}. 
Amer. J. Math. (to appear).\par

[DGP] W. Decker, G.-M. Greuel, G. Pfister, {\psl Primary decomposition:
algorithms and comparisons}, in {\pit Algorithmic Algebra and Number Theory}, 187-220,
Springer Verlag, Heidelberg, 1998.\par  

[Ei] D. Eisenbud.  {\pit Commutative Algebra with a View Toward Algebraic
Geometry}, Springer-Verlag, 1995.\par  

[EH] D. Eisenbud, J. Harris. {\pit The geometry of schemes},
Springer-Verlag, 1999.\par   

[EHV] D. Eisenbud, C. Huneke, W. Vasconcelos, {\psl Direct methods for primary
decomposition}. Invent. Math. {\pbo 110} (1992), 207--235.\par

[Gr] M. Green, {\psl Generic initial ideals},
in {\pit Six lectures on commutative algebra},
 Birkh{\"a}user, Progress in Mathematics {\pbo 166},
(1998), 119--185.\par

[GLP] L. Gruson, R. Lazarsfeld, C. Peskine, {\psl On a theorem of Castelnuovo,
and the equations defining space curves}. Invent. Math. {\pbo 72} (1983),
491--506.\par 

[Hu] H. Hulett, {\psl Maximum Betti numbers of homogeneous ideals with a
given Hilbert
function}.  Comm. in Algebra {\pbo 21} (1993), 2335--2350.\par

[KP]  T. Krick, L. M. Pardo, {\psl A computational method for Diophantine
approximation}. Algorithms in algebraic geometry and applications
(Santander, 1994), 193--253, {\pit Progr. Math.}, {\pbo 143}, Birkh{\"a}user,
Basel, 1996.\par   

[L] R. Lazarsfeld, {\psl Multiplier ideals for algebraic geometers}.
Notes for the ICTP School on Vanishing Theorems, August
2000.\par
 
[MV] C. Miyazaki, W. Vogel, {\psl
Bounds on cohomology and Castelnuovo-Mumford regularity}.
J. Algebra {\pbo 185} (1996), 626--642.\par

[Mu] D. Mumford, {\pit Lectures on curves on an algebraic surface}, 
Ann. of Math. Stud. {\pbo 59}, Princeton University Press (1966).
\par

[NS1] U. Nagel, P. Schenzel, {\psl Cohomological annihilators and 
Castelnuovo-Mumford regularity}.  Contemp. Math. {\pbo 159} (1994),
307--328.\par

[NS2] U. Nagel, P. Schenzel, {\psl Degree bounds for generators of 
cohomology modules and Castel\-nuovo-Mumford regularity}. Nagoya
Math. J. {\pbo 152} (1998), 153--174.\par

[Pa] K. Pardue, {\psl
Deformation classes of graded modules and maximal Betti numbers}.
Illinois J. Math. {\pbo 40} (1996), 564--585. \par

[V1] W. Vasconcelos, {\psl Constructions in commutative algebra,} in {\pit
  Computational algebraic geometry and commutative algebra  (Cortona,
  1991)}, 151--197, Sympos. Math., XXXIV, Cambridge Univ. Press,
  Cambridge, 1993.\par

[V2] W. Vasconcelos, {\pit Arithmetic of blowup algebras.} London
Mathematical Society Lecture Note Series, 195. Cambridge 
University Press, Cambridge, 1994.\par

[V2] W. Vasconcelos, {\pit Computational methods in commutative algebra and
algebraic geometry}.  Algorithms and Computation in Mathematics, 2.
Springer-Verlag, Berlin, 1998.\par 

}
\bigskip\bigskip\bigskip

\noindent Marc Chardin, Institut de Math{\'e}matiques,
CNRS \&\ Universit{\'e} Paris 6,

4, place Jussieu, F--75252 Paris {\nrm cedex} 05, France

chardin@math.jussieu.fr

\end

[PS] C. Peskine, L. Szpiro, {\psl Liaison des vari{\'e}t{\'e}s
  alg{\'e}briques. I}. Invent. Math. {\pbo 26} (1974), 271-302.\par

[SY] T. Shimoyama, K. Yokoyama. Localization and primary decomposition of
polynomial ideals. {\psl J. Symbolic Comput.} {\pbo 22} (1996), no. 3, 247--277. 
\par

[GTZ] P. Gianni, B. Trager, G. Zacharias. Gr{\"o}bner bases and primary
decomposition of polynomial ideals. Computational aspects of commutative
algebra. {\psl J. Symbolic Comput.} {\pbo 6} (1988), no. 2-3, 149--167.
\par

A very remarkable and strange result follows from the work of Nadel,
Ein, Lazarsfeld and others on multiplier ideal (this is implicitely in
[L], as it grew out of a discussion with Dale Cutkosky): 
\medskip
{\bf Theorem 3.4.} {\sl If $k$ is of characteristic zero, there exists an
  ideal $J\subseteq I^{*}$ such that 
  $\sqrt {I}=\sqrt{J}$ and 
$$
\reg (J)\leq n(d_{1}-1)+1.
$$
}

It seems that one may impose $J$ to sit in $I$ and get the 
bound $\reg (J)\leq (n+1)(d_{1}-1)+1$. In any case, this is extremely
surprising when one remembers that all results and conjectures up to
now where are proving or expecting good bounds to occur under nice
geometric hypotheses on the support : This result implies that if you
allow to have embedded 
components at the singular locus, you can have  bounds close to the
expected ones on regularity... If you are (very) optimistic and
believe on the conjecture: $\reg (\sqrt{J})\leq \reg (J)$, then
it implies that $\reg  (\sqrt{I})\leq n(d_{1}-1)+1$ with no hypothesis
on the geometry! This in turn would imply that the radical of an ideal
is computable in small degree (the same type of bound, replacing $n$
by a bigger constant in $n$), using an algorithm based on the ideas of
[EVH, V1]. 

Also, a refined version of this result implies most of the previous
bounds on regularity, and some improvements. However, the arguments
relies on sophisticated  techniques that are unkonwn or fail in
positive characteristic (Hironaka's version of resolution of
singularities and refined versions of Kodaira's vanishing theorem, in
particular). 
\medskip
